\documentclass[11pt]{article}
\usepackage{geometry}                
\usepackage[parfill]{parskip}    
\usepackage{amssymb}
\usepackage{epstopdf}
\usepackage{amsmath,amsthm,amscd,amssymb}
\usepackage{latexsym}
\numberwithin{equation}{section}

\theoremstyle{plain}
\newtheorem{theorem}{Theorem}[section]
\newtheorem{lemma}[theorem]{Lemma}
\newtheorem{corollary}[theorem]{Corollary}

\newtheorem{conjecture}[theorem]{Conjecture}

\theoremstyle{definition}
\newtheorem{definition}[theorem]{Definition}

\theoremstyle{remark}
\newtheorem{remark}[theorem]{Remark}

\newtheorem{case[theorem]}{Case}

\newcommand {\binomial}[2]{\left( \begin{array}{c} #1 \\#2 \end{array}
\right)}

\title{Erd\"os-Falconer distance problem, exponential sums, and Fourier analytic approach to incidence theorems in vector spaces over finite fields}
\author{Alex Iosevich and Doowon Koh}

\begin{document}

\maketitle

\begin{abstract} We study the Erd\"os/Falconer distance problem in  vector spaces over finite fields with respect to the cubic metric. Estimates for discrete Airy sums and Adolphson/Sperber estimates for  
exponential sums in terms of Newton polyhedra play a crucial role.  Similar techniques are used to study the incidence problem between points and cubic and quadratic curves. As a result we obtain a non-trivial range of exponents that appear to be difficult to attain using combinatorial methods. \end{abstract}

\tableofcontents

\section{Introduction}

\subsection{The Erd\"os distance problem} The Erdos distance  
conjecture in the Euclidean space says that if $E$ is a finite subset  
of ${\Bbb R}^d$, $d \ge 2$, then
\begin{equation} \label{erdosconjecture}  \# \Delta(E) \gtrapprox  
{(\# E)}^{\frac{2}{d}}, \end{equation} where
$$ \Delta(E)=\{|x-y|: x,y \in E\}, $$ with ${|x-y|}^2={(x_1-y_1)}^2+ 
\dots+{(x_d-y_d)}^2$ and here, and throughout the paper, $X \lesssim Y 
$ means that there exists $C>0$ such that $X \leq CY$, and $X  
\lessapprox Y$, with the controlling parameter $N$, means that for  
every $\epsilon>0$ there exists $C_{\epsilon}>0$ such that $X \leq C_ 
{\epsilon}N^{\epsilon}Y$.

Taking $E={\Bbb Z}^d \cap {[0,N^{\frac{1}{d}}]}^d$ shows that (\ref 
{erdosconjecture}) cannot in general be improved. The conjecture has  
not been solved in any dimension. See, for example, (\cite{Ma02}),  
(\cite{AP95}), and the references contained therein for the description  
of the conjecture, background material, and a survey of recent results.

In this paper we study the Erd\"os distance problem in vector spaces  
over finite fields. This problem was
recently addressed by Tao (\cite{T}) who relates it to some  
interesting questions in combinatorics, and, more recently, by  
Iosevich and Rudnev. We shall describe these results later in the  
introduction.

Let ${\mathbb F}_q$ denote the finite field with $q$ elements, and let ${\mathbb F}^d_q$ denote the 
$d$-dimensional vector space over this field. Let $E \subset {\mathbb F}_q^d$, $d \ge 2$. Then the  
analog of the classical Erd\"os distance problem is to determine the  
smallest possible cardinality of the set
$$ \Delta_n(E)=\{{||x-y||_n=(x_1-y_1)}^n+\dots+{(x_d-y_d)}^n: x,y \in E 
\}, $$ with $n$ a positive integer
$\ge 2$, viewed as a subset of ${\mathbb F}_q$.

In the finite field setting, the estimate (\ref{erdosconjecture})  
cannot hold without further restrictions.
To see this, let $E={\mathbb F}^d_q$. Then $\# E=q^d$ and $\# \Delta 
(E)=q$. Furthermore, an interesting feature of the Erd\"os distance  
problem in the finite field setting with $n=2$ is the existence of  
non-trivial spheres of $0$ radius. These are sets of the form $\{x  
\in {\mathbb F}_q^d: x_1^2+x_2^2+\dots+x_d^2=0\}$ and several  
assumptions in the statements of results below are there precisely to  
deal with issues created by the presence of this object. For example,  
suppose $-1$ is a square in ${\mathbb F}_q$. Using spheres of radius  
$0$ one can show, in even dimensions, that there exists a set of  
cardinality precisely $q^{\frac{d}{2}}$ such that all the distances, $ 
{(x_1-y_1)}^2+\dots+{(x_d-y_d)}^2$ are $0$. What's more, suppose $ 
{\mathbb F}_q$ is a finite field, such that $q=p^2$, where $p$ is a  
prime. Then $E={\mathbb F}_p^d$ is naturally embedded in ${\mathbb F} 
_q^d,$ has cardinality $q^{\frac{d}{2}}$, and determines only $\sqrt 
{q}$ distances. If $n>2$, the situation is equally fascinating. For  
example, if $n=3$ and $d=2$, the equation $x_1^3+x_2^3=0$ always has  
at least $q$ solutions, since cube root of $-1$ is $-1$. This  
equation may have as many as $3q$ solutions if the primitive cube  
root of $-1$ is in the field.

With these examples as guide, we generalize the conjecture originally  
stated in (\cite{IR06}) in the case $n=2$ as follows.

\begin{conjecture} \label{conjecture} Let $E \subset {\mathbb F}^d_q$  
of cardinality $\ge  Cq^{\frac{d}{2}}$, with $C$ sufficiently large.  
Then
$$ \# \Delta_n(E) \gtrapprox q. $$ \end{conjecture}

The authors conjecture in (\cite{IR06}) that the constant $C$ that   appears above may be taken to be any number bigger than one, at least in the case $n=2$. It is   interesting to note that if $n>2$, the situation becomes more   complicated. For example, as we pointed out above, if $n=3$ and $d=2 $, the number of points on the curve $x_1^3+x_2^3=0$ may be as high as $3q$, depending on whether or not the primitive cube root of $-1$  is in the field. Thus a corresponding conjecture in the case $n>2$  
must be designed with these issues in mind.

\section{Previous results} A Euclidean plane argument due to Erd\"os  
(\cite{E45}) can be applied to the finite field set-up under the  
assumption of Conjecture \ref{conjecture} to show that if $d=2$ and $ 
\# E \ge Cq$, with $C$ sufficiently large, then

\begin{equation} \label{trivial2d} \# \Delta_n(E) \gtrsim {(\# E)}^ 
{\frac{1}{2}}. \end{equation}

This result was improved by Bourgain, Katz and Tao (\cite{BKT04}) who  
showed using intricate incidence geometry that for every $\epsilon>0, 
$ there exists $\delta>0,$ such that if $\# E \lesssim q^{2-\epsilon} 
$, then
$$ \# \Delta_2(E) \gtrsim q^{\frac{1}{2}+\delta}. $$

The relationship between $\epsilon$ and $\delta$ in the above argument is difficult to determine. Moreover, matters are even more subtle in higher dimensions in the context of vector spaces over finite fields because intersection  of analogs of spheres, both quadratic and cubic, in ${\mathbb F}^d_q$ may be quite complicated,  and the standard induction on the dimension argument in ${\mathbb R}^d$ 
(see e.g. \cite{AP95}) that allows one to bootstrap the estimate (\ref {trivial2d}) into the estimate
\begin{equation} \label{trivialhighd} \# \Delta_{{\mathbb R}^d}(E)  \gtrsim {(\# E)}^{\frac{1}{d}} \end{equation} does not immediately go  through. We establish the finite field analog of the estimate (\ref 
{trivialhighd}) below using Fourier analytic methods and number  theoretic properties of Kloosterman sums and its more general analogs. 

Another way of thinking of Conjecture \ref{conjecture} is in terms of  the Falconer distance conjecture
(\cite{Fal85}) in the Euclidean setting which says that if the  Hausdorff dimension of a set in ${\mathbb R}^d$ exceeds $\frac{d}{2} $, then the Lebesgue measure of the distance set is positive. Conjecture \ref{conjecture} implies that if the size of the set is  greater than $q^{\frac{d}{2}}$, then the distance set contains a positive proportion of all the possible distances, an analogous statement.

In (\cite{IR06}) the authors proved the following result.

\vskip.125in 

\begin{theorem} \label{iosevichrudnev} Let $E \subset {\Bbb F}_q^d$,  
$d \ge 2$, such that $\# E \ge Cq^{\frac{d+1}{2}}$. Then if $C$ is  
sufficiently large, $\Delta_2(E)$ contains every element of 
${\Bbb F} 
_q$. \end{theorem}

\vskip.125in 

\section{Main results of this paper}

\subsection{Distances determined by a single set}

Our first result is the version of Theorem \ref{iosevichrudnev} for  
cubic metrics.

\vskip.125in 

\begin{theorem} \label{main} Suppose that $q$ is a prime number  
congruent to $1$ modulo $3$. Let $E \subset {\Bbb F}_q^d$, such that $ 
\# E \ge Cq^{\frac{d+1}{2}}$. Then if $C$ is sufficiently large, $ 
\Delta_3(E)$ contains every element of ${\Bbb F}_q$. 

Suppose that $d=2$, and $n \ge 2$. Then if $\# E \ge Cq^{\frac{3}{2}}$ for $C$ sufficiently large, then 
$\Delta_n(E)$ contains every elements of ${\Bbb F}_q$.

\end{theorem}

\vskip.125in 

\begin{corollary} \label{score} Suppose that $q$ is a prime number  
congruent to $1$ modulo $3$. Let $E \subset {\Bbb F}_q^d$, $d \ge 2$,  
such that $\# E=Cq^{\frac{d+1}{2}}$. Then if $C$ is sufficiently large,
$$ \# \Delta_3(E) \approx {(\# E)}^{\frac{2}{d+1}}.$$

In two dimensions, the same conclusions, with $d=2$, holds for any $n \ge 2$. 
\end{corollary}

Note that in the case $d=2$, the exponent $\frac{2}{3}$ obtained via the corollary, for the given range of parameters, is a much better  exponent than the one obtained by the incidence argument due to Erdos  
described in (\ref{trivial2d}) above. Also, we point out once more that Erdos' argument does not generalize to higher dimensions, at least not very easily, due to the possibly complicated intersection properties of cubic varieties. 

\subsection{Szemeredi-Trotter type Incidence theorems and distances between pairs of sets} As  
in the case $n=2$, the proof of Theorem \ref{main} can be modified to yield a good upper bound on the  
number of incidences between points and cubic surfaces in vector spaces over finite fields. It is an analog, and a higher dimensional generalization, of the following classical result due to Szemeredi and Trotter. 

\vskip.125in

\begin{theorem} \label{st} The number of incidences between $N$  
points and $M$ lines (or circles of the same radius) in the plane is
$$ \lesssim N+M+{(NM)}^{\frac{2}{3}}. $$
\end{theorem}

\vskip.125in 

Our incident estimate is the following. 

\vskip.125in

\begin{theorem} \label{incidence3} Suppose that $q$ is a prime number congruent to 1 modulo 3.
Let $E,F \subset {\Bbb F}_q^d$, $d\ge 2$. Then if  $j \not=0$,
$$ \# \{(x,y) \in E \times F: {(x_1-y_1)}^3+\dots+{(x_d-y_d)}^3=j \} $$
$$ \lesssim \# E \cdot \# F \cdot q^{-1}+q^{\frac{d-1}{2}} \cdot {(\# E)} 
^{\frac{1}{2}} \cdot {(\# F)}^{\frac{1}{2}}. $$

Similarly, if $q$ is a prime number and $j \not=0$, then
$$ \# \{(x,y) \in E \times F: {(x_1-y_1)}^2+\dots+{(x_d-y_d)}^2=j \} $$
$$ \lesssim \# E \cdot \# F \cdot q^{-1}+q^{\frac{d-1}{2}} \cdot {(\# E)} 
^{\frac{1}{2}} \cdot {(\# F)}^{\frac{1}{2}}. $$

In two dimensions, the same result holds, with $d=2$, with $\Delta_3$ replaced by $\Delta_n$ for any 
$n \ge 2$. 

\end{theorem}

\begin{remark} In particular, if $\# E \approx \# F \approx q^{\frac{d+1}{2}}$, then the number of incidences between points in $E$ and  "spheres", quadratic or cubic, centered at elements of $F$ is $\lesssim q^d$.

To make the numerology more transparent, Theorem \ref{incidence3}  says that if $N \approx q^{\frac{d+1}{2}}$, the number of incidences  between $\approx N$ points and $\approx N$ spheres, cubic or quadratic, in ${\Bbb F}_q^d$  is $\lesssim q^d=N^{\frac{2d}{d+1}}$. In two dimensions this says  that the number of incidences between $N$ points and $N$ circles is $ \lesssim N^{\frac{4}{3}}$, provided that $N \approx q^{\frac{d+1}{2}}$, matching in this setting the exponent in the celebrated result due to Szemeredi and Trotter in the Euclidean plane (see Theorem \ref{st} above). 
\end{remark} 

\vskip.125in 

An easy modification of the method used to prove Theorem \ref{incidence3} above yields the following distance set result. 

\vskip.125in 

\begin{corollary} \label{twosets} Let $E,F \subset {\Bbb F}_q^d$, $d  
\ge 2$. Suppose that $q$ is a prime number congruent to 1 modulo 3 and 
$\# E \cdot \# F \ge Cq^{d+1}$. Let $\Delta_3 (E,F)=\{||x-y||_3: x \in E, y \in F\}$.
Then if $C$ is sufficiently large, then
$\Delta_3(E, F)$ contains every element of ${\Bbb F}^{*}_q$.

As before, in two dimensions the same conclusion holds, with $d=2$, with $\Delta_3$ replaced by 
$\Delta_n(E)$. 
\end{corollary}

Observe that if $E=F$, then we can safely say that in fact $\Delta_3(E, F)$ contains every element of 
${\Bbb F}_q$, but if $E \not=F$, the zero distance may not be present. 

We also call the reader's attention to the fact that an analogous version of this result was independently obtained by Shparlinski in \cite{Sh06}. 

\vskip.125in

\section{Fourier analytic preliminaries and notation}

\vskip.125in

Let ${\Bbb F}_q$ be a finite field with $q$ elements, where $q$ is a  
prime number. Let
$$ \chi(t)=e^{\frac{2 \pi i}{q}t}. $$

Given a complex valued function $f$ on ${\Bbb F}_q^d$, define the  
Fourier transform of $f$ by the equation
$$ \widehat{f}(m)=q^{-d} \sum_{x \in {\Bbb F}_q^d} \chi(-x \cdot m) f 
(x). $$

We also need the following basic identity, typically known as the Plancherel theorem. Let $f$ be as above. Then
$$ \sum_{m \in {\Bbb F}_q^d} {|\widehat{f}(m)|}^2=q^{-d} \sum_{x \in  
{\Bbb F}_q^d} {|f(x)|}^2.$$

\vskip.125in

\section{Proof of the first part of Theorem \ref{main}}

\vskip.125in

Let $\chi(s)=e^{\frac{2 \pi i}{q}s}$. Let $S_j$ denote the  
characteristic function of the cubic "sphere"
$$ \{x \in {\Bbb F}_q^d: ||x||_3=j\}, $$ where, as above,
$$ ||x||_3=x_1^3+ \dots+ x_d^3.$$

The key estimate of the paper is the following. 

\vskip.125in 

\begin{theorem} \label{cubic} Let $||x||_3=x_1^3+ \dots+x_d^3 $. 
Suppose that $q$ is a prime number congruent to 1 modulo 3 and $j \not =0$. Then if  
$m \not=(0, \dots, 0)$, then
$$ \left| \widehat{S}_j(m) \right| =\left| q^{-d} \sum_{\{x \in {\Bbb F}_q^d: ||x||_3=j\}} \chi(x \cdot  
m) \right| \lesssim q^{-\frac{d+1}{2}}, $$ and if $m =(0, \dots, 0)$, then
$$ \widehat{S}_j (m) = q^{-1} + O(q^{-\frac{d+1}{2}})$$
$$ \approx  q^{-1} .$$
\end{theorem}

For $j\ne 0$, consider
$$ \# \{(x,y) \in E \times E: ||x-y||_3=j\}$$
$$=\sum_{x,y\in {\Bbb F}_q^d} E(x)E(y)S_j(x-y)$$
$$=q^{2d} \sum_m {|\widehat{E}(m)|}^2 \widehat{S}_j(m)=I+II, $$ where
$$ I=q^{2d} {|\widehat{E}(0, \dots, 0)|}^2 \widehat{S}_j(0, \dots, 0), $$ and
$$II=q^{2d} \sum_{m \not=(0, \dots, 0)} {|\widehat{E}(m)|}^2 \widehat{S}_j(m).$$

Using the second part of Theorem \ref{cubic},
$$ I\approx q^{2d} q^{-2d} {(\# E)}^2 \cdot q^{-1} .$$

Whereas using the first part of Theorem \ref{cubic},
$$ |II| \lesssim q^{2d} q^{-\frac{d+1}{2}} \sum_{m \not=(0,\cdots, 0)}
{|\widehat{E}(m)|}^2$$
$$ \lesssim q^{2d} q^{-\frac{d+1}{2}} q^{-d} \sum_{x \in {\Bbb F}_q^d}E^2(x)$$
$$=q^{\frac{d-1}{2}} \cdot \# E.$$

We therefore obtain that
$$  \# \{(x,y) \in E \times E: ||x-y||_3=j\}=I+II, $$ where
$$ I \gtrsim {(\# E)}^2 q^{-1}, $$ and
$$ |II| \lesssim \# E \cdot q^{\frac{d-1}{2}}. $$

We conclude that if $\# E \ge Cq^{\frac{d+1}{2}}$, with $C$  
sufficiently large, then
$$  \# \{(x,y) \in E \times E: ||x-y||_3=j\}>0$$ for each $j \not=0$.  
This completes the proof of Theorem \ref{main}.

\section{Proof of Theorem \ref{cubic}}

We have
$$ \widehat{S}_j(m)=q^{-d} \sum_{\{x \in {\Bbb F}_q^d: ||x||_3=j\}} \chi 
(-x \cdot m)$$
$$=q^{-1}\delta(m)+q^{-d-1} \sum_x \sum_{t \in {\Bbb F}_q^ 
{*}} \chi(t(||x||_3-j)) \chi(-x \cdot m), $$ where $\delta(m)=1$ if $m= 
(0, \dots, 0)$ and $0$ otherwise. 

\vskip.125in 

\begin{lemma} \label{cubicn0} Let $\chi$ be a nontrivial additive character of $F_q$ with $q \equiv 1 \mod(3).$  Suppose that $m=(m_1,\cdots, m_l) \in {({\Bbb F}_q^{*})}^l.$ Then for any multiplicative
character $\psi$ of $F_q$ of order 3 and $t \not=0,$ we have
$$ \prod_{j=1}^{l} \sum_{s_j \in F_q} \chi(-s_j m_j +s_j^3 t) $$
$$=\psi^{-l}(t) \sum_{s_1,\cdots, s_l \in F_q^*}
 \chi(s_1+\cdots+s_l+m_1^3t^{-1}s_1^{-1}+\cdots+m_l^3t^{-1}s_l^{-1} ) \psi(s_1)\cdots \psi(s_l) ,$$ 
where $3^{-3}m_j^{3}$ is denoted by $m_j^{3}$ in the right-hand side of the equation.\end{lemma}

\vskip.125in 

We shall also need the following result due to Duke and Iwaniec (\cite {DI93}).

\vskip.125in 

\begin{theorem} \label{dukeiwaniec} Suppose that $q \equiv 1 \mod(3)$ and let $\psi$ be a multiplicative character of order three. Then
$$ \sum_{s \in {\Bbb F}_q} \chi(as^3+s)=\sum_{s \in {\Bbb F}_q^{*}}  \psi(sa^{-1}) \chi(s-{(3^3as)}^{-1}), $$ for any $a \in {\Bbb F}_q^{*}$. \end{theorem}
It follows that
$$ \sum_{s \in {\Bbb F}_q} \chi(-sm_j+s^3t)=\sum_{s \in {\Bbb F}_q}  \chi(s-s^3tm_j^{-3})$$
$$=\sum_{s \in {\Bbb F}_q^{*}} \psi(st^{-1}) \chi(s+m_j^3 t^{-1}3^{-3} 
s^{-1}).$$ since $\psi$ is a multiplicative character of $F_q$ of order three and $m_j\not=0.$ Absorbing $3^{-3}$ into $m_j$ to make the notations simple, we complete the proof of Lemma \ref{cubicn0}.

\vskip.125in

\begin{lemma} \label{cubic0} Let $\chi$ be a nontrivial additive character of $F_q$ with $q \equiv 1 \mod(3).$ Then for any multiplicative character $\psi$ of $F_q$ of order 3 and $t \not=0,$ we have 
$$ {\Big( \sum_{s\in F_q} \chi(ts^3) \Big)}^l = \sum_{r=0}^l \binomial{l}{r} q^l \,\psi^{-(l+r)}(t) 
{\Big(\widehat{\psi}(-1) \Big)}^{l-r} \,{\Big(\widehat{\psi^2}(-1) \Big)}^{r} ,$$
where $\binomial{.}{.} $ is a binomial coefficient , $l$ is a positive integer,
and the Fourier transform of a multiplicative characer $\psi$ of ${\Bbb F}_q$ is given by
$$ \widehat{\psi}(v)= q^{-1}\sum_{s\in {\Bbb F}_q^*}\chi(-vs)\psi(s).$$ \end{lemma}

\begin{remark} $\widehat{\psi}(v) =O(q^{-\frac{1}{2}})$ for $v \ne 0$.
\end{remark}

To prove Lemma \ref{cubic0} , we need the following theorem.
For the proof,  see the (\cite{RH}, page 217, Theorem 5.30).

\vskip.125in 

\begin{theorem} \label{cubice} Let $\chi$ be a nontrivial additive character of $F_q ,\,n\in \mathbb{N} ,$ and $\psi$ a multiplicative character of $F_q$ of order $ h=$gcd$(n,q-1).$ Then
$$\sum_{s\in F_q} \chi(ts^n+b)  = \chi(b) \sum_{k=1}^{h-1} \psi^{-k}(t)\, G(\psi^k, \chi)$$
for any $t,b\in F_q$ with $t\not=0$, where $G(\psi^k, \chi) =\sum_{s\in F_q^*} \psi^k(s) \chi(s).$
\end{theorem}

By using Theorem \ref{cubice}, we see that for any multiplicative character $\psi$ of order three,
$$ {\Big(\sum_{s\in F_q} \chi(ts^3) \Big)}^l $$
$$={\Big(\sum_{k=1}^{2} \psi^{-k}(t)\, \sum_{s\in F_q^*} \psi^k(s)\, \chi(s) \Big)}^l$$
$$ ={\Big(\psi^{-1}(t) \,\sum_{s\in F_q^*} \psi(s)\,\chi(s)+ \psi^{-2}(t)
\sum_{s\in F_q^*} \psi^2(s) \chi (s) \Big)}^l$$
$$= {\Big( G_1(t)+G_2(t) \Big)}^l$$
$$ =\sum_{r=0}^l \binomial{l}{r} {G_1(t)}^{l-r} {G_2(t)}^r,$$
where 
$$G_1(t)= \psi^{-1}(t) \sum_{s\in F_q^*} \psi(s) \chi(s)$$ and 
$$G_2(t)= \psi^{-2}(t) \sum_{s\in F_q^*} \psi^2(s)\chi (s).$$ 

Note that $G_1(t) = q \psi^{-1}(t) \,\widehat{\psi}(-1)$ and $ G_2(t) = q\psi^{-2}(t)\, \widehat{{\psi}^2}(-1).$ 

Thus we conclude that
$${\Big(\sum_{s\in F_q} \chi(ts^3) \Big)}^l =\sum_{r=0}^l \binomial{l}{r} q^l \psi^{-(l+r)}(t) 
{\Big(\widehat{\psi}(-1) \Big)}^{l-r} {\Big(\widehat{\psi^2}(-1) \Big)}^{r}.$$

We are now ready to prove Theorem \ref{cubic}. First, we assume that $m=(0,\cdots,0) \in F_q^d$. Then, using Lemma \ref{cubic0}, we see that
$$ \widehat{S}_j(0,\cdots,0)=q^{-d} \sum_{\{x \in {\Bbb F}_q^d: ||x||_3=j\}} 1$$
$$=q^{-1}+q^{-d-1} \sum_{t \in {\Bbb F}_q^{*}} \chi(-tj)\,\sum_x  \chi(t(||x||_3))$$ 
$$= q^{-1}+ q^{-d-1}\sum_{t \in {\Bbb F}_q^{*}} \chi(-tj)\, \sum_{r=0}^d \binomial{d}{r} q^d \psi^{-(d+r)}(t) 
{\Big(\widehat{\psi}(-1) \Big)}^{d-r} {\Big(\widehat{\psi^2}(-1) \Big)}^{r} $$
$$= q^{-1}+ q^{-1}\sum_{r=0}^d \binomial{d}{r} {\Big(\widehat{\psi}(-1) \Big)}^{d-r} \,{\Big(\widehat{\psi^2}(-1) \Big)}^{r} \sum_{t \in {\Bbb F}_q^{*}} \chi(-tj) \psi^{-(d+r)}(t)$$
$$= q^{-1}+ q^{-1}\sum_{r=0}^d \binomial{d}{r} {\Big(\widehat{\psi}(-1) \Big)}^{d-r} \,{\Big(\widehat{\psi^2}(-1) \Big)}^{r} q \widehat{\psi^{-(d+r)}}(j)$$
$$ =q^{-1}+ O(q^{-\frac{d+1}{2}}) \approx  q^{-1} .$$ 

\vskip.125in 

In the last equality, we used the fact that $ \widehat{\psi}(v) = O(q^{-\frac{1}{2}})$
for any multiplicative character of $F_q$ with $v\ne 0$. 
Thus the second part of Theorem \ref{cubic} is proved.

In order to prove the first part of Theorem \ref{cubic},
we shall deal with the problem in case $ m=(m_1,\cdots,m_d)
\not = (0,\cdots,0)$. Suppose that $m_j \not=0 $ for $j\in J \subset \{1,2,\cdots,d\} $ 
and $m_j=0 $ for $j \in \{1,2,\cdots,d\}\setminus J =J^{'}.$
Without loss of generality, we may assume that 
$J=\{1,2,\cdots,l\}$ and $J^{'} = \{l+1,\cdots, d\}$ for some $l=1,2,\cdots,d.$
Using Lemma \ref{cubicn0} and Lemma \ref{cubic0}, we see that
$$\widehat{S_j}(m) = q^{-d-1}\sum_{t\in F_q^*} \chi(-tj)
\sum_{x\in F_q^d} \chi(t||x||_3 - m\cdot x) $$
$$=  q^{-d-1}\sum_{t\in F_q^*} \chi(-tj)
\Big( \prod_{k=1}^l\sum_{s_k \in F_q}\chi(ts_k^3 - m_k s_k) \Big)
\Big(\prod_{k=l+1}^d \sum_{s_k\in F_q}\chi(ts_k^3)\Big)$$
$$=q^{-d-1}\sum_{t\in F_q^*} \chi(-tj)\psi^{-l}(t) \sum_{s_1,\cdots, s_l \in F_q^*}
 \chi(s_1+\cdots+s_l+m_1^3t^{-1}s_1^{-1}+\cdots+m_l^3t^{-1}s_l^{-1} ) \psi(s_1)\cdots \psi(s_l)  $$
$$ \times \sum_{r=0}^{d-l} \binomial{d-l}{r} q^{d-l}
\psi^{-(d-l+r)}(t)
{\Big(\widehat{\psi}(-1) \Big)}^{d-l-r} {\Big(\widehat{\psi^2}(-1) \Big)}^{r} $$
$$=q^{-1-l}\sum_{r=0}^{d-l}\binomial{d-l}{r}{\Big(\widehat{\psi}(-1)\Big)}^{d-l-r}
 {\Big( \widehat{\psi^2}(-1) \Big)}^r \sum_{t\in F_q^*} \chi(-tj) \psi^{-(d+r)}(t)$$
 $$ \times \sum_{s_1,\cdots, s_l \in F_q^*}
 \chi(s_1+\cdots+s_l+m_1^3t^{-1}s_1^{-1}+\cdots+m_l^3t^{-1}s_l^{-1} ) \psi(s_1)\cdots \psi(s_l). $$
Since $\binomial{d-l}{r}{\Big(\widehat{\psi}(-1)\Big)}^{d-l-r}
{\Big( \widehat{\psi^2}(-1) \Big)}^r  = O(q^{-\frac{1}{2}(d-l)})$, we obtain that
$$ \Big|\widehat{S_j}(m) \Big| \lesssim q^{-1-\frac{d+l}{2}} \sum_{r=0}^{d-l}\,|A_r(\chi,\psi) |,$$
where $A_r(\chi,\psi)$ is given by 
$$\sum_{t \in F_q^*} \chi(-tj)\psi^{-(d+r)}(t) \sum_{s_1,\cdots, s_l \in F_q^*}
 \chi(s_1+\cdots+s_l+m_1^3t^{-1}s_1^{-1}+\cdots+m_l^3t^{-1}s_l^{-1} ) \psi(s_1)\cdots \psi(s_l) .$$

We now apply the result of Adolphson and Sperber  
(\cite{AS89}, Theorem 4.2, Corollary 4.3) to see that for all $r=0,1,\cdots, d-l,$
$$ |A_r(\chi,\psi)| \lesssim q^{\frac{l+1}{2}} .$$

This completes the proof. 
                                                
\vskip.125in

\section{Proof of the second part of Theorem \ref{main}} 

As in the proof of the first part of Theorm \ref{main}, it suffices to prove the follwoing estimation.

\begin{theorem} \label{2D} Let $||x||_n=x_1^n+ x_2^n $
for $x \in \Bbb{F}_q^2$ and $n\ge2$.
Suppose that $q$ is a prime number and $j \not =0$. Then if
$m \not=(0,0) $, then
$$ \left| \widehat{S}_j(m) \right| =\left| q^{-2} \sum_{\{x \in {\Bbb F}_q^2: ||x||_n=j\}} \chi(-x \cdot  
m) \right| \lesssim q^{-\frac{3}{2}}, $$
and if $m =(0, 0)$,then
$$ \widehat{S}_j (m) = q^{-1} + O(q^{-\frac{3}{2}}) \approx  q^{-1} .$$
\end{theorem}

To prove Theorem \ref{2D}, we observe that for $j\ne 0$ and $m \in \Bbb{F}_q^2$,

$$ \widehat{S}_j(m)=q^{-2} \sum_{\{x \in {\Bbb F}_q^2: ||x||_n=j\}} \chi(-x \cdot m)$$
$$=q^{-1}\delta(m)+q^{-3} \sum_x \sum_{t \in {\Bbb F}_q^ 
{*}} \chi(t(||x||_n-j)) \chi(-x \cdot m), $$ where $\delta(m)=1$ if $m= 
(0,  0)$ and $0$ otherwise.

First we shall prove the second part of Theorem \ref{2D}.
Using Theorem \ref{cubice}, we see that for a multi-index
$\beta =(\beta_1,\cdots,\beta_{h-1})$,
$$\Big(\sum_{s\in\Bbb{F}_q}\chi(ts^n)\Big)^2$$
$$=\sum_{\beta_1+\cdots+\beta_{h-1}=2} \frac{2!}{\beta_1!\cdots \beta_{h-1}!}
\psi^{-(\beta_1+\cdots+(h-1)\beta_{h-1})}(t)q^2 {\Big(\widehat{\psi}(-1)\Big)}^{\beta_1}
\cdots {\Big(\widehat{\psi^{h-1}}(-1)\Big)}^{\beta_{h-1}} $$
where $\psi$ is a multiplicative character of $\Bbb{F}_q$ of order
$h= gcd(n, q-1)$.
It therefore follows that
$$\widehat{S}_j(0,0)
=q^{-1}+\sum_{\beta_1+\cdots+\beta_{h-1}=2}\frac{2!}{\beta_1!\cdots \beta_{h-1}!}
\widehat{\psi^{-\gamma(h,\beta)}}(j)
{\Big(\widehat{\psi}(-1)\Big)}^{\beta_1}\cdots
{\Big(\widehat{\psi^{h-1}}(-1)\Big)}^{\beta_{h-1}}$$
where $\gamma(h,\beta)$ is given by $\beta_1+2\beta_2+\cdots+(h-1)\beta_{h-1}.$

Since $\widehat{\psi}(v)= O(q^{-\frac{1}{2}})$
for each multiplicative character $\psi$ and $v\in \Bbb{F}_q^*$,
we conclude 
$$\widehat{S}(0,0) = q^{-1}+O(q^{-\frac{3}{2}})\approx q^{-1}.$$
This completes the proof of the second part of Theorem \ref{2D}.

It remains to prove the first part of Theorem \ref{2D}.
The cohomological interpretation can be used to estimate the exponential sums.
We now introduce the cohomology theory based on work of authors
in \cite{JF} and \cite{AS89}.
Let $g$ be a polynomial given by
\begin{equation}\label{Polynomial} g = \sum_{\alpha\in J}A_{\alpha} x^{\alpha}
\in \Bbb{F}_q [x_1,\cdots,x_d] ,\end{equation}
where $J$ is a finite subset of $(\Bbb{N} \cup \{0\})^d$,
and $A_{\alpha}\ne 0$ if $ \alpha \in J.$
We denote by $\sum(g)$ the Newton polyhedron of $g$ which is the convex hull
in $\Bbb{R}^d$ of the set $J\cup (0,\cdots,0)$.
For any face $\sigma$ (of any dimension) of $\sum(g)$, we put
$$ g_{\sigma} = \sum_{\alpha \in \sigma \cap J} A_{\alpha}x^{\alpha}.$$

\begin{definition}Let $g\in \Bbb{F}_q [x_1,\cdots,x_d]$ be a polynomial as in (\ref{Polynomial}).
We say that $g$ is nondegenerate with respect to $\sum(g)$ if for every face $\sigma$
of $\sum(g)$ that does not contain the origin, the polynomials 
$$ \frac{\partial{g_{\sigma}}}{\partial{x_1}} , \cdots,
\frac{\partial{g_{\sigma}}}{\partial{x_d}}$$
have no common zero in $\Big({\bar{\Bbb{F}_q}^*}\Big)^d$
where $\bar{\Bbb{F}_q}$ denotes an algebraic clousure of $\Bbb{F}_q.$
We say that $g$ is commode with respect to $\sum(g)$ if for each
$k=1,2,\cdots,d$, $g$ contains a term $A_kx_k^{\alpha_k}$ for some
$\alpha_k > 0$ and $ A_k \ne 0.$\end{definition}

The general version of the following theorem can be found in \cite{JF}
(see Theorem 9.2).

\vskip.125in 

\begin{theorem}\label{Cohomology} Let $q$ be a prime number.
Suppose that $g : \Bbb{F}_q^d \rightarrow \Bbb{F}_q,
d\ge 2$, is commode and nondegenerate with respect to $\sum(g)$. Then
$$ \sum_{x\in \Bbb{F}_q^d} \chi(g(x)) = O(q^{\frac{d}{2}}).$$
\end{theorem}

We now prove the first part of Theorem \ref{2D}.
Since $m \ne (0,0)$, we have $ \sum_{x\in \Bbb{F}_q^2} \chi(-x\cdot m) =0$.
We therefore see that for $j\ne 0$,
$$\widehat{S_j}(m)
= q^{-3}\sum_{(t,x_1,x_2) \in \Bbb{F}_q^*\times \Bbb{F}_q^2}\chi(g(t,x_1,x_2))
= q^{-3}\sum_{(t,x_1,x_2) \in \Bbb{F}_q^3}\chi(g(t,x_1,x_2)),$$
where $g(t,x_1,x_2) =tx_1^n+tx_2^n-m_1x_1-m_2x_2-jt.$

If $m_1\cdot m_2 \ne 0$, then $g$ is commode. By Theorem \ref{Cohomology},
it suffices to show that $g$ is nondegenerate with respect to $\sum(g)$.
Note that $\sum(g)$ has five zero-dimensional faces, eight one-dimensional
faces and three two-dimensional faces which do not contain the origin.
It is easy to show that for every face $\sigma$ of $\sum(g)$ that does not
contain the origin, the polynomials 
$$ \frac{\partial{g_{\sigma}}}{\partial{t}}, \frac{\partial{g_{\sigma}}}{\partial{x_1}},
\frac{\partial{g_{\sigma}}}{\partial{x_2}}$$
have no common zero in ${(\Bbb{F}_q^*)}^3$ because we may assume that 
$q$ is sufficiently large and so $n$ is not congruent to 0 modulo $q$.
This implies that $g$ is nondegenerate with respect to $\sum(g)$.
We now assume that $ m_1 \cdot m_2 =0$.
Without loss of generality, we may assume that $m_1\ne 0$, and
$m_2=0$ because $m\ne (0,0).$
By using Theorem \ref{cubice}, we obtain that for a multiplicative
character $\psi$ of $\Bbb{F}_q$ of order $h = gcd (n, q-1)$,
$$\widehat{S_j}(m)
=q^{-3}\sum_{(t,x_1)\in \Bbb{F}_q^*\times \Bbb{F}_q}\chi(tx_1^n-m_1x_1-jt)
\sum_{k=1}^{h-1}\psi^{-k}(t) q \widehat{\psi^k}(-1)$$
$$=q^{-2}\sum_{k=1}^{h-1}\widehat{\psi^k}(-1)
\sum_{(t,x_1)\in \Bbb{F}_q^*\times \Bbb{F}_q}\psi^{-k}(t)\chi(tx_1^n-m_1x_1-jt)$$
$$\lesssim q^{-2}q^{-\frac{1}{2}}\sum_{k=1}^{h-1}|R_k(\psi^{-k},\chi)|,$$
where $R_k(\psi^{-k},\chi)$ is given by
$$\sum_{(t,x_1)\in \Bbb{F}_q^*\times \Bbb{F}_q}\psi^{-k}(t)\chi(tx_1^n-m_1x_1-jt).$$
For each $k=1,2,\cdots,h-1$, define $\psi^{-k}(0)=0$. Then we can obtain that
$$R_k(\psi^{-k},\chi)
=\sum_{(t,x_1)\in \Bbb{F}_q\times \Bbb{F}_q}\psi^{-k}(t)\chi(tx_1^n-m_1x_1-jt).$$
Applying Theorem \ref{Cohomology}, we have
$$R_k(\psi^{-k},\chi)=O(q).$$
This completes the proof.

\vskip.125in 

\section{Proof of Theorem \ref{incidence3} and Corollary \ref{twosets}}

As we mentioned in the introduction, this is a simple variation on the
proof of Theorem \ref{main}. Indeed,
$$ \# \{(x,y) \in E \times F: ||x-y||_n=j\}$$
$$=q^{2d} \sum_m \overline{\widehat{E}(m)} \widehat{F}(m) \widehat{S} 
_j(m)$$
$$=\# E \cdot \# F \cdot \widehat{S} 
_j(0,\cdots,0) +q^{2d} \sum_{m \not=(0,  
\dots, 0)} \overline{\widehat{E}(m)} \widehat{F}(m) \widehat{S}_j(m)=I 
+II.$$

By the second part of Theorem \ref{cubic}(or Theorem \ref{2D}),
$$ I \lesssim \# E \cdot \# F \cdot q^{-1}. $$

Applying Cauchy-Schwartz, Theorem \ref{cubic}( or Theorem \ref{2D}) and Plancherel, we see  
that
$$ |II| \lesssim q^{2d} q^{-\frac{d+1}{2}} \sum_{m \not=(0, \dots, 0)}
|\widehat{E}(m)| |\widehat{F}(m)| $$
$$ \leq q^{2d} q^{-\frac{d+1}{2}} {\left( \sum_m {|\widehat{E}(m)|}^2  
\right)}^{\frac{1}{2}} \cdot {\left( \sum_m {|\widehat{F}(m)|}^2  
\right)}^{\frac{1}{2}}$$
$$ \leq q^{2d} q^{-\frac{d+1}{2}} q^{-d} {\left( \sum_x {|E(x)|}^2  
\right)}^{\frac{1}{2}} \cdot {\left( \sum_x {|F(x)|}^2  
\right)}^{\frac{1}{2}}$$
$$=q^{\frac{d-1}{2}} \cdot \sqrt{\# E} \cdot \sqrt{\# F}.$$
This completes the proof of Theorem \ref{incidence3}.
\vskip.125in
In order to prove Corollary \ref{twosets}, we observe that by
the second part of Theorem \ref{cubic}(or Theorem \ref{2D}),
$$ I \gtrsim \# E \cdot \# F \cdot q^{-1}.$$
On the other hand, we have seen above that
$$ |II| \lesssim q^{\frac{d-1}{2}} \cdot \sqrt{\# E} \cdot \sqrt{\# F},  
$$ and the result follows by a direct comparison.


\vskip.125in 




\newpage

\begin{thebibliography}{17}

\vskip.125in

\bibitem{AP95} J. Pach, and P. Agarwal, {\it Combinatorial geometry,}  
Wiley-Interscience Series in Discrete Mathematics and Optimization. A  
Wiley-Interscience Publication. John Wiley and Sons, Inc., New York  
(1995).

\bibitem{AS89} A. Adolphson and S. Sperber, {\it Exponential sums and  
Newton polyhedra: cohomology and estimates}, (1989), \textbf{130},  
367-406.

\bibitem{BKT04} J. Bourgain, N. Katz, and T. Tao, {\it A sum-product  
estimate in finite fields, and applications}, Geom. Funct. Anal.  
\textbf{14} (2004), 27-57.

\bibitem{DI93} W. Duke and H. Iwaniec, {\it A relation between cubic  
exponential and Kloosterman sums}, Contemp. Math. {\bf 143}, (1993),  
255-258.

\bibitem{E45} P. Erd\"os {\it On sets of distances of n points,}  
Amer. Math. Monthly. \textbf{53} (1946), 248--250.

\bibitem{Fal85} K. J. Falconer. {\it On the Hausdorff dimensions of  
distance sets.} Mathematika {\bf 32} (1985), 206-212.

\bibitem{G03} B. J. Green, {\it Restriction and Kakeya phonomena,}  
Lecture notes (2003).

\bibitem{IK04} H. Iwaniec, and E. Kowalski, {\it Analytc Number  
Theory,} Colloquium Publications \textbf{53} (2004).

\bibitem{IR06} A. Iosevich and M. Rudnev, {\it Erdos/Falconer  
distance problem in vector spaces over finite fields}, TAMS (to  
appear), (2006).

\bibitem{Katz88} N. Katz, {\it Gauss sums, Kloosterman sums, and  
monodromy groups,} Ann. Math. Studies {\textbf 116},
Princeton (1988).

\bibitem{LM06} M. Lacey and W. McClain, {it On an argument of Shkredov in the finite field setting}, (2006), On-line journal of analytic combinatorics (http://www.ojac.org). 

\bibitem{Ma02} J. Matousek, {\it Lectures on Discrete Geometry,}  
Graduate Texts in Mathematics, Springer \textbf{202} (2002).

\bibitem{MT04} G. Mockenhaupt, and T. Tao,
{\it Restriction and Kakeya phenomena for finite fields,} Duke Math.  
J. \textbf{121} (2004), 35--74.

\bibitem{N91} H. Niederreiter, {\it The distribution of values of  
Kloosterman sums,} Arch. Math. \textbf{56} (1991), 270--277.

\bibitem{Sh06} I. Shparlinski, {\it On the set of distances between  
two sets in vector spaces over finite fields}, (2006), (preprint).

\bibitem{StSh03} E. Stein, and R. Shakarchi, {\it Fourier analysis,}  
Princeton Lectures in Analysis, (2003).

\bibitem{T} T. Tao, {\it Finite field analogues of Erd\"os, Falconer,  
and Furstenberg problems,} preprint.

\bibitem{We48} A. Weil, {\it On some exponential sums,} Proc. Nat.  
Acad. Sci. U.S.A. \textbf{34}
(1948), 204--207.

\bibitem{RH} R. Lidl and H. Niederreiter, {\it Finite fields},
 Cambridge Univ. Press (1997).

\bibitem{JF} J. Denef, and F. Loeser,{\it Weights of exponential sums,
intersection cohomology, and Newton polyhedra,} Invent. Math.
\textbf{106} (1991), 275--294.


\end{thebibliography}
\end{document}